\documentclass[letterpaper,    %Option: Papier
               english,        %Option: Sprache
               12pt,           %Option: Schriftgröße
               twoside]        %Option: beidseitig
               {amsart}        %Paket:  (American Mathematical Society) Stil des Dokuments
\usepackage[english]{babel}                   %Paket:  Trennungsregeln
\usepackage[T1]                               %Option: T1-Fonts verwenden
           {fontenc}                          %Paket:  Fonts
\usepackage{ae}                               %Paket:  schmalere CM-T1-Fonts verwenden
\usepackage[letterpaper,                      %Option: Papiergröße
            bookmarks,                        %Option: erstellt Lesezeichen
            bookmarksopen=true,               %Option: öffnet Lesezeichen beim Start
            bookmarksnumbered=true,           %Option: nummeriert Lesezeichen
            pdfauthor={David Klawonn and
                       Christof Puhle},       %Option: Autor des Dokuments
            pdftitle={Remarks on ''Resolving
                      isospectral 'drums'
                      by counting nodal
                      domains''},             %Option: Titel des Dokuments
            pdfkeywords={flat 4-tori,
                         isospectral tori,
                         nodal domains,
                         Laplace operator},   %Option: Schlüsselwörter
            colorlinks,                       %Option: verwendet farbige Verweise statt Kästchen
            linkcolor=black,                  %Option: Farbe der internen Verweise
            citecolor=black,                  %Option: Farbe der Zitate
            urlcolor=black,                   %Option: Farbe der URLs
            plainpages=false]                 %Option: verhindert spezielle Fehlermeldungen
            {hyperref}                        %Paket:  Hyperlinks/Lesezeichen
\usepackage[nobysame]                         %Option: schaltet "by same author"-Linie aus
           {amsrefs}                          %Paket:  (American Mathematical Society) Literaturverzeichnis            
\usepackage{amstext}                          %Paket:  (American Mathematical Society) angepasster Text in Gleichungen
\usepackage{amssymb}                          %Paket:  (American Mathematical Society) zusätzliche Symbole
\advance\oddsidemargin by -0.93cm
\advance\evensidemargin by -0.92cm
\advance\topmargin by -0.55cm
\theoremstyle{plain}
\newtheorem*{thm}{Theorem}   
\begin{document}
\thispagestyle{empty}
\date{September 24, 2007}
\title[Remarks on ``Resolving isospectral `drums' by counting ...'']{Remarks on
``Resolving isospectral `drums' by counting nodal domains''}
\author{Jochen Bruening, David Klawonn and Christof Puhle}
\address{\hspace{-5.5mm}
{\normalfont\ttfamily bruening@mathematik.hu-berlin.de}\newline
{\normalfont\ttfamily klawonn@mathematik.hu-berlin.de}\newline
{\normalfont\ttfamily puhle@mathematik.hu-berlin.de}\newline
Institut f\"ur Mathematik \newline
Humboldt-Universit\"at zu Berlin\newline
Unter den Linden 6\newline
10099 Berlin, Germany}
\thanks{Supported by GIF and the SFB 647: `Space--Time--Matter'}
\subjclass[2000]{Primary 58 J 53; Secondary 58 J 50}
\begin{abstract}
In \cite{GnSm05} the authors studied the $4$-parameter family of
isospectral flat $4$-tori $T^\pm(a,b,c,d)$ discovered by Conway and
Sloane. With a particular method of counting nodal domains they were
able to distinguish these tori (numerically) by computing the corresponding
nodal sequences relative to a few explicit tuples $(a,b,c,d)$. In this note
we confirm the expectation expressed in \cite{GnSm05} by proving analytically
that their nodal count distinguishes any $4$-tuple of distinct positive real
numbers.
\end{abstract}
\maketitle
\pagestyle{headings}
%
%
%
%------------------------------------------------------------------------------
%
\section{Introduction}\noindent
In 1964 J.~Milnor \cite{Mil64} constructed two $16$-dimensional non-isometric
flat tori with the same spectrum for the Laplace-Beltrami operator on forms
of every degree, and thus produced the first example of non-isometric {\em isospectral}
manifolds. Since then many examples of such manifolds (see for example \cite{Sun85},
\cite{BCDS94} and \cite{CoSl92}) have been found and studied.

While it remains still unclear to what extent the spectrum of the Laplace-Beltrami
operator determines the geometry of the underlying manifold, these examples show
that the spectrum does not contain enough information to determine the manifold
and its metric uniquely. It has been proposed recently that the {\em nodal count},
i.e.~the number of nodal domains of the eigenfunctions of the Laplace-Beltrami
operator might provide the missing information, such that spectrum and nodal count
together should yield isometry. Indeed, in \cite{GnSm05} the authors used the
$4$-parameter family of isospectral flat $4$-tori $T^\pm(a,b,c,d)$ constructed by
Conway and Sloane \cite{CoSl92} to show how the isospectrality can be `resolved'
using nodal domains. By counting the latter in very special way, and then arranging
the result in a so-called `nodal sequence' they were able to exhibit that these nodal
sequences for the tori belonging to four carefully chosen tuples $(a,b,c,d)$ are
different, if numerically.

In this note we shall give an alternative and analytic way to show that the nodal
sequence defined in \cite{GnSm05} distinguishes every pair of isospectral tori in
this family if all four parameters are distinct; if at least two of them are equal
then the corresponding tori are isometric (cf.~\cite{CoSl92}*{p.~94, Remark 2}).
\begin{thm}
For any choice of distinct positive numbers $a,b,c,d\in\mathbb{R}_{+}$
the nodal sequences of the tori $T^+(a,b,c,d)$ and $T^-(a,b,c,d)$
are distinct.
\end{thm}
The paper is structured as follows: In \autoref{sec:2} we collect
some facts on flat $n$-tori, their spectra and their eigenfunctions, and
define the notions of `nodal domain' and `nodal sequence'. We then proceed
to introduce the isospectral flat tori of Conway and Sloane in \autoref{sec:3},
and eventually prove the theorem stated above.
%
%
%
%------------------------------------------------------------------------------
%
\section{Nodal sequences of flat tori}\label{sec:2}\noindent
Let $v_1,\dots,v_n\in\mathbb{R}^n$ denote linearly independent vectors
and
\begin{equation*}
\Gamma:=\mathrm{span}_{\mathbb{Z}}\left\{v_1,\dots,v_n\right\}
\end{equation*}
the lattice generated by these vectors. The flat torus given by the
lattice $\Gamma$ is
\begin{equation*}
T:=\mathbb{R}^n/A\mathbb{Z}^n
\end{equation*}
where the columns of the $(n\times n)$-matrix $A$ consist of the vectors $v_i$
\begin{equation*}
A=[v_1,\dots,v_n].
\end{equation*}
The \emph{Gram matrix} of $T$ is defined by $G:=A^{\top}A$ and $Q:=G^{-1}$ denotes
its inverse. The regular matrix $Q$ determines the torus completely, $T=T(Q)$.
The dual lattice is
\begin{equation*}
\Gamma^*:=\mathrm{span}_{\mathbb{Z}}\left\{v^*_1,\dots,v^*_n\right\}
\end{equation*}
with $(v^*_i)_i$ the dual basis,  $v^*_i(v_j)=\delta_{ij}$.

The \emph{Laplace-Beltrami operator} $\Delta$ on $T$ takes the form
\begin{equation*}
\Delta=-\sum_{i=1}^n\frac{\partial^2}{\partial x_i^2}.
\end{equation*}
Its spectrum relative to the torus $T(Q)$ consists only of isolated
eigenvalues with finite multiplicity and can be computed explicitly,
\begin{equation*}
\mathrm{spec}_{T}(\Delta)=\left\{4\,\pi^2\,q^{\top}\,Q\,q\, : \,
q\in\mathbb{Z}^n\right\}.
\end{equation*}
A certain eigenvalue $\lambda\in\mathrm{spec}_{T}(\Delta)$ may correspond
to multiple \emph{representing vectors}, i.e.~vectors $q\in\mathbb{Z}^n$
satisfying $\lambda=4\,\pi^2\,q^{\top}\,Q\,q$. The number of distinct
representing vectors relative to the eigenvalue $\lambda$ is called
the \emph{degeneracy} of $\lambda$. The degeneracy of a given
$\lambda$ equals the dimension of its eigenspace, a basis of which
is given by the functions
\begin{equation*}
\Psi_{q}:T \ni x\longmapsto\exp\left( 2\,\pi\, i\,
\sum_{i=1}^nq_i\,v^*_i(x)\right)\in\mathbb{C}
\end{equation*}
where $q=(q_1,\dots,q_n)^{\top}\in\mathbb{Z}^n$ is a representing vector
of $\lambda$.

Let $f:M\rightarrow\mathbb{R}$ be a function on a compact manifold $M$.
Then the \emph{nodal domains} of $f$ are defined as the connected components
of $M\setminus f^{-1}(0)$, the number of which is finite.

Throughout this work we will only consider the nodal domains of the real and
imaginary parts of the eigenfunctions $\Psi_{q}$. We will count these
domains in the same way as introduced in \cite{GnSm05}. Here is the procedure:
First split $\Psi_{q}$ into its real and imaginary parts
\begin{equation*}
\Psi_{q}^\mathrm{re}(x)=\cos\left(2\,\pi\,\sum_{i=1}^nq_i\,v^*_i(x)\right),
\quad\quad
\Psi_{q}^\mathrm{im}(x)=\sin\left(2\,\pi\,\sum_{i=1}^nq_i\,v^*_i(x)\right).
\end{equation*}
Then introduce the transformation
\begin{equation*}
T(Q)\ni x\longmapsto Q^{-1}\,x=y\in \tilde{T}
\end{equation*}
onto the standard torus $\tilde{T}=\mathbb{R}^n/\mathbb{Z}^n$, which
gives rise to the Laplace-Beltrami operator
\begin{equation*}
\tilde{\Delta}=-\sum_{i,j}Q_{ij}\,\frac{\partial^2}{\partial y_i\partial y_j}
\end{equation*}
and to the functions
\begin{equation*}
\tilde{\Psi}_{q}^\mathrm{re}(y)=\cos\left(2\,\pi\,q^{\top}\,y\right), \quad\quad
\tilde{\Psi}_{q}^\mathrm{im}(y)=\sin\left(2\,\pi\,q^{\top}\,y\right).
\end{equation*}
The number of nodal domains is -- by definition -- the number given by
lifting these functions to $\mathbb{R}^n$ and then counting their nodal
domains in the unit cube ignoring identifications at the boundary.
The resulting number which we call the \emph{nodal count} $\nu(q)$ for a
given representing vector $q\in\mathbb{Z}^n$ is given by the following
formula (see \cite{GnSm05}):
\begin{equation*}
\nu(q)=
\left\{\begin{array}{lcl}
	 2\,\sum_{i=1}^n|q_i| & \mathrm{for} & \Psi_{q}^\mathrm{im}\\[1mm]
	 2\,\sum_{i=1}^n|q_i| + 1 & \mathrm{for} & \Psi_{q}^\mathrm{re}.
\end{array} \right.
\end{equation*}

Since $T(Q)$ is a compact manifold we can arrange its
spectrum $\mathrm{spec}_{T}(\Delta)$ in increasing order
\begin{equation*}
0<\lambda_1<\lambda_2<\dots<\lambda_i<\dots.
\end{equation*}
If we compute the nodal count of every vector $q\in\mathbb{Z}^n$,
a finite set of nodal counts $\{\nu^i_1,\nu^i_2,\dots\}$ belongs to
each eigenvalue $\lambda_i$. The cardinality of this set equals
the degeneracy of the corresponding eigenvalue. One obtains the
\emph{nodal sequence}
\begin{equation*}
\left\{\left\{\nu^1_1,\nu^1_2,\dots\right\},\left\{\nu^2_1,\nu^2_2,\dots\right\},
\dots,\left\{\nu^i_1,\nu^i_2,\dots \right\},\dots\right\}
\end{equation*}
by fitting each nodal sequence in the same position as the corresponding
eigenvalue in the spectrum. By means of this sequence we shall
distinguish isospectral tori.
%
%
%
%------------------------------------------------------------------------------
%
\section{The construction of Conway and Sloane}\label{sec:3}\noindent
Our work deals with the $4$-parameter family of isospectral flat
tori $T^\pm(a,b,c,d)$ discovered by Conway and Sloane \cite{CoSl92}.
As mentioned in the previous section, these tori are described
by the inverse $Q^\pm(a,b,c,d)$ of the corresponding Gram matrix.
Explicitly,
\begin{equation*}
Q^+=\frac{1}{12}\left[\hspace{-1mm}
\begin{array}{cccc}
       9a+b+c+d  &  3a-3b-c+d &  3a+b-3c-d&  3a-b+c-3d \\
       3a-3b-c+d &  a+9b+c+d  &  a-3b+3c-d&  a+3b-c-3d \\ 
       3a+b-3c-d &  a-3b+3c-d &  a+b+9c+d &  a-b-3c+3d \\
       3a-b+c-3d &  a+3b-c-3d &  a-b-3c+3d&  a+b+c+9d  \\
\end{array}\hspace{-1mm}\right]\hspace{-1mm},
\end{equation*}
\begin{equation*}
Q^-=U^{\top}Q^+U \quad\mathrm{with}\quad
U=\frac{1}{2}\left[\begin{array}{cccc}
-1&1&1&1\\
-1&-1&-1&1\\
-1&1&-1&-1\\
-1&-1&1&-1\end{array}\right].
\end{equation*}
The defining parameters $a,b,c,d$ are required to be strictly positive.
It is remarked in \cite{CoSl92} that the tori $T^+$ and $T^-$ are equivalent
if two of these parameters are equal. Therefore we shall only consider vectors
$(a,b,c,d)$ of pairwise distinct positive numbers. We are now ready for the
\begin{proof}[Proof of the Theorem]
To begin with, we define for $m\in\mathbb{N}$ the set 
\begin{equation*}
V_{m}:=\left\{(q_1,q_2,q_3,q_4)^{T}\in\mathbb{Z}^4\, : \,\sum_i|q_i|=m\right\}.
\end{equation*}
This set is obviously finite and contains all vectors $q\in\mathbb{Z}^4$
that represent the nodal count $2m$ or $2m+1$, according to whether we
consider $\Psi_q^\mathrm{im}$ or $\Psi_q^\mathrm{re}$. We then define
\begin{equation*}
E_{m}^{\pm}:=\left\{4\,\pi^2\, q^\top\, Q^{\pm}\, q\,:\,q\in V_m\right\}
\end{equation*}
as the set of eigenvalues with a representing vector of nodal count $2m$ or $2m+1$.
Thus, if $E_{m}^{+}$ and $E_{m}^{-}$ do not coincide for a certain $m$, then
the nodal sequences of the tori $T^+$ and $T^-$ are distinct.

The $E_{m}^{\pm}$'s can be viewed as sets of linear functions in the variables
$a,b,c,d$. By inspection we obtain equality ($E_m^+=E_m^-$) for $m=1,2,3$. The
first interesting case appears for $m=4$, where
\begin{eqnarray*}
E^{+}_{4} & =& (4\,\pi^2/3)
\Big\{ (4 a + 25 b + c),(25a + b + 4c),(a + 4 b + 25 c),\\
&& (b + 25 c + 4 d),(25 a + 4 b + d),(25 b + 4c + d),\\
&& (4 a + 25c + d),(4 a + b + 25 d),(a + 25 b + 4 d), \\
&& (25 a + c + 4 d),(4 b + c + 25d),(a + 4 c + 25d),\\
&& (4 a + 16 b + 9 c + d),(9 a + 4 b + 16 c + d),\\
&& (16 a + 9 b + 4 c + d),(9 a + 16 b + c + 4 d),\\
&& (16 a + b + 9c + 4 d),(16 a + 4 b + c + 9 d),\\
&& (a + 16 b + 4 c + 9 d),(4 a + b + 16 c + 9 d),\\
&& (4 a + 9 b + c + 16 d),(a + 4 b + 9 c + 16 d),\\
&& (a + 9 b + 16 c +4 d),(9 a + b + 4 c + 16 d) \Big\} \,\, \cup \,\, (E^{+}_{4}\cap E^{-}_{4}),
\end{eqnarray*}
\begin{eqnarray*}
E^{-}_{4}& =& (4\,\pi^2/3)
\Big\{ (25 a + 4 b + c),(a + 25 b + 4 c),(4 a + b + 25 c),\\
&& (4 a + 25 b + d),(25 a + 4 c + d),(4 b + 25 c + d),\\
&& (25 a + b + 4 d),(25 b + c + 4 d),(a + 25 c + 4 d),\\
&& (a + 4 b + 25 d),(4 a + c + 25 d),(b + 4 c + 25 d),\\
&& (9 a + 16 b + 4c + d),(16 a + 4b + 9 c + d),\\
&& (4 a + 9 b + 16 c + d),(16 a + 9 b + c + 4 d),\\
&& (a + 16 b + 9 c + 4 d),(4 a + 16 b + c + 9 d),\\
&& (16 a + b + 4 c + 9 d),(a + 4 b + 16 c + 9 d),\\
&& (9 a + 4 b + c + 16 d),(4 a + b + 9 c + 16 d),\\
&& (9 a + b + 16 c + 4 d),(a + 9 b + 4 c + 16 d) \Big\} \,\, \cup \,\,
(E^{+}_{4}\cap E^{-}_{4}).
\end{eqnarray*}
Inspecting the sets $E_{4}^{+}$ and $E^{-}_{4}$ more carefully one
notes that there are two sets of coefficients -- namely $(1,4,9,16)$
and $(0,1,4,25)$ -- such that $E^{+}_{4}$ contains all even permutations
of the variables $a,b,c,d$ in the linear forms with these coefficients
while $E^{-}_{4}$ contains the odd ones. Hence we may assume that $a<b<c<d$
and obtain a unique maximum among all elements of $E^{+}_{4}\cup E^{-}_{4}$,
namely $b+4c+25d$. Consequently $E^{+}_{4} \neq E^{-}_{4}$ as claimed.
\end{proof}
%
%
%
%------------------------------------------------------------------------------
%
\section*{Acknowledgement}\noindent
We want to thank Uzy Smilansky for inspiring discussions and both the German-Israeli
Foundation GIF and the SFB 647: \emph{Space--Time--Matter} for financial support.
%
%
%
%------------------------------------------------------------------------------
%

%
%
%
\end{document}